\documentclass{amsart}
\usepackage{amssymb,amscd}
\advance\hoffset by -1in
\advance\voffset by -1in

\usepackage[cp1251]{inputenc}
\usepackage[english,russian]{babel}
\newtheorem{theorem}{Theorem}
\newtheorem{lemma}{Lemma}

\newtheorem{proposition}{Proposition}
{\par\noindent{\bf Proof.}} 
{\hfill$\scriptstyle\blacksquare$}
\title{The theorems on freedom for relatively free Lie algebras with a relations}

\author{A.\,F.\,Krasnikov}
\address{Omsk State University, pr. Mira 55-A, 644077, Omsk, Russia}
\email{phomsk@mail.ru}

\begin{document}
\selectlanguage{english}
\maketitle
\markboth{A.\,F.\,Krasnikov}{The theorems on freedom}

\section*{Introduction}
If $M$ is an ideal of the algebra $L$, generated
by the set $X$, then we write $M=\mbox{\rm id}_L (X)$.
The well-known result of Shirshov \cite{Sh} tells us that if
$F$ --- a free Lie algebra on free generators $y_1,\ldots,y_n$,
$H$ --- the subalgebra of $F$, generated by $y_1,\ldots,y_{n-1}$, $r\in F\setminus H$, $R=\mbox{\rm id}_F (r)$
then $H\cap R=0$.

Talapov \cite{Tl} proved an analogous result for polynilpotent Lie algebras.

In \S1 is contained the definition of Fox's derivations and various considerations in the universal enveloping algebra of the free Lie algebra. The prove the theorem on freedom for relatively free Lie algebras with a single relation occupies \S2:
\begin{theorem}\label{tm2_alg}
Suppose $F$ is a free Lie algebra on free generators $y_1,\ldots,y_n\,(n>2)$,
$H$ --- the subalgebra of $F$, generated by $y_1,\ldots,y_{n-1}$, $N_{11}$  ---  an ideal of $F$, $F/N_{11}$  ---  a relatively free Lie algebra,
\begin{eqnarray}\label{vved}
N_{11} \geqslant \ldots \geqslant N_{1,m_1+1}=N_{21} \geqslant \ldots \geqslant N_{s,m_s+1},
\end{eqnarray}
where $N_{kl}$ --- the $l$-th term of the lower central series of $N_{k1}$.
Let $r$ be an element of $N_{1i}\backslash N_{1,i+1}\,(i\leqslant m_1)$, $R=\mbox{\rm id}_F (r)$.
If (and only if) $r\not\in H+N_{1,i+1}$ then $H\cap (R+N_{kl})=H\cap N_{kl}$ for each term $N_{kl}$ of series {\rm (\ref{vved})}.
\end{theorem}

From the theorem \ref{tm2_alg} we have the theorem on freedom for polynilpotent Lie algebras with a single relation \cite{Tl}.

Kharlampovich \cite{Hm} proved that if $F$ --- a free Lie algebra on free generators $y_1,\ldots,y_n$,
$r_1,\ldots,r_m$ --- an elements of $F$ $(m<n)$, $R=\mbox{\rm id}_F (r_1,\ldots,r_m)$,
then there exists a subset $y_{j_1},\ldots,y_{j_p}$ $(p\geqslant n-m)$ of $y_1,\ldots,y_n$
such that $H\cap R=0$, where $H$ is the subalgebra of $F$ generated by $y_{j_1},\ldots,y_{j_p}$.

Kharlampovich \cite{Hm} proved also an analogous result for solvable Lie algebras.

The prove the theorem on freedom for relatively free Lie algebras with a relations (generalized Freiheitssatz) occupies \S3:
\begin{theorem}\label{tm3_alg}
Suppose $F$ is a free Lie algebra on free generators $y_1,\ldots,y_n$,
$N_{11}$  --- an ideal of $F$, $F/N_{11}$ --- solvable Lie algebra,
\begin{eqnarray}\label{end_algr_3}
N_{11} \geqslant \ldots \geqslant N_{1,m_1+1}=N_{21} \geqslant \ldots \geqslant N_{s,m_s+1},
\end{eqnarray}
where $N_{kl}$  --- the $l$-th term of the lower central series of $N_{k1}$.
Let $r_1,\ldots,r_m$ be an elements of $N_{11}$ $(m<n)$, $R=\mbox{\rm id}_F (r_1,\ldots,r_m)$. Then there exists a subset $y_{j_1},\ldots,y_{j_p}$ $(p\geqslant n-m)$ of $y_1,\ldots,y_n$
such that $H\cap (R+N_{kl}) = H\cap N_{kl}$ for each term $N_{kl}$ of series {\rm (\ref{end_algr_3})}, where $H$ is the subalgebra of $F$ generated by $y_{j_1},\ldots,y_{j_p}$.
\end{theorem}

From the theorem \ref{tm3_alg} we have the theorem on freedom for polynilpotent Lie algebras with a relations.


\section{Some properties of Fox's derivations}
We regard all the algebras under consideration as given
over an arbitrarily fixed field $P$. Let $L$ be a Lie algebra, $u$, $v$ be an elements of $L$.
We denote by $[u,v]$ the multiplication on $L$ and by $L_{(k)}$ the $k$-th term of the lower central series of $L$.

We denote by $U(L)$ the universal enveloping algebra of $L$
and by $U_0(L)$ --- the ideal in $U(L)$ which is generated by $L$.
If $M$ is an ideal in $L$ then we denote by $M_U$ the ideal in $U(L)$ which is
generated by $M$. The algebra $L$ is identified with a suitable subalgebra in $U(L)^{(-)}$ and
\begin{eqnarray}\label{sob_pr}
uv=vu + [u,v].
\end{eqnarray}

Let $A$, $B$ be a sets of elements of $U(L)$. Then we define the subset $AB$
of $U(L)$ by $AB=\{ab \mid a\in A,\,b\in B \}$.

It is worthy of note that a subalgebra of a free Lie algebra is a free Lie algebra itself \cite{Sh2}.

The Poincar\'e--Birkhoff--Witt theorem tells us that if $u_1,\ldots,u_n,\ldots$ is an ordered basis in $L$ then $U(L)$ has for
basis $1$ and the set of monomials of the form $u_{i_1}\cdots u_{i_r}$ where $i_1\leqslant\ldots\leqslant i_r$.

Let $F$ be a free Lie algebra with a free set $\{g_j \mid  j\in J\}$ of generators.
If $u\in U(F)$ then we can find unique elements $D_j(u)\in U(F)$ such that
\begin{eqnarray*}
u = \sum_{j\in J} g_jD_j(u).
\end{eqnarray*}
We call the elements $D_j(u)$ $(j\in J)$ the Fox derivatives of $u$.
It is easy to deduce the following relations:
\begin{gather}
D_k(\alpha u+\beta v)=\alpha D_k(u)+\beta D_k(v)\,(k\in J);\notag\\
D_j(g_j)=1\,(j\in J),\,D_k(g_j)=0,\mbox{ if }k\neq j;\notag\\
D_k([u,v])= D_k(u)v-D_k(v)u;\notag\\
D_k([n,u])\equiv D_k(n)u\mod{N_U}\,(k\in J),\label{four}
\end{gather}
where $N$ --- an ideal of $F$; $u,\,v\in F;\,n\in N;\,\alpha,\,\beta\in P$.

Let $H$ be a subalgebra of $F$ with a free set $\{h_k \mid  k\in K\}$ of generators, $\{\partial_k \mid  k\in K\}$ ---  the Fox derivatives of $U(H)$, $f\in H$.\\
From $f=\sum_{k\in K} h_k\partial_k(f)$ and $h_k=\sum_{j\in J} g_jD_j(h_k)$ follows that\\
$f=\sum_{j\in J} g_j(\sum_{k\in K} D_j(h_k)\partial_k(f))$. Therefore
\begin{eqnarray}\label{alg1_dcf}
D_j(f)=\sum_{k\in K} D_j(h_k)\partial_k(f),\,j\in J.
\end{eqnarray}

\begin{lemma}\label{alg1_lm1}
Suppose $F$ is a free Lie algebra on free generators $\{g_j \mid  j\in J\}$, $N$ --- an ideal of $F$,
$D_j\,(j\in J)$ --- the Fox derivatives of $U(F)$, $K\subseteq J$, $F_K$
the subalgebra of $F$, generated by $\{g_j\mid  j\in K\}$, $\{u_j\mid j\in K\}$ an
elements of $U(F_K)$ such that $u_j$ is equal to zero for all but
a finite number of $j$. If
\begin{eqnarray}\label{alg1_lm1_1}
u=\sum_{j\in K} g_ju_j\equiv 0\mod{N_U},
\end{eqnarray}
then there is an element $v\in F_K\cap N$ such that $D_j(v)\equiv u_j \mod{N_U}$ $(j\in K)$.
\end{lemma}
\begin{proof}
Let $\{a_i,\,i\in J_1\}$ be an ordered basis in $F_K\cap N$,
$\{a_i,\,i\in J_1\}\cup\{b_j,\,j\in J_2\}$ an ordered basis in
$F_K$, $\{a_i,\,i\in J_1\}\cup\{b_j,\,j\in J_2\}\cup\{c_s,\,s\in
J_3\}$ an ordered basis in $F_K+N$ ($\{c_s,\,s\in J_3\}\subseteq N$) and $\{a_i,\,i\in
J_1\}\cup\{b_j,\,j\in J_2\}\cup\{c_s,\,s\in J_3\}\cup\{d_t,\,t\in
J_4\}$ an ordered basis in $F$. We assume that $a_i<b_j<c_s<d_t$
for any $i\in J_1,\,j\in J_2,\,s \in J_3,\,t\in J_4$. Then
$U(F)$ has for basis the set of monomials of the form
\begin{eqnarray}\label{alg1_1}
a_{i_1}\ldots a_{i_\mu}b_{j_1}\ldots b_{j_\nu}c_{s_1}\ldots
c_{s_\eta}d_{t_1}\ldots d_{t_\theta},
\end{eqnarray}
where $i_1\leqslant\ldots\leqslant i_\mu,\,j_1\leqslant\ldots\leqslant
j_\nu,\,s_1\leqslant\ldots\leqslant s_\eta,\,t_1\leqslant\ldots\leqslant
t_\theta,\,\mu\geqslant 0,\,\nu\geqslant 0,\,\eta\geqslant 0,\,\theta\geqslant 0$.
From $u\in U(F_K)$ follows that $u$ is a linear combination of monomials of the form (\ref{alg1_1})
with $\eta=\theta=0$. From $u\in N_U$ follows that $u$ is a linear combination of monomials of the form (\ref{alg1_1})
with $\mu+\eta\geqslant 1$. Whence $u$ is a linear combination of monomials of the form (\ref{alg1_1})
with $\mu\geqslant 1$, $\eta=\theta=0$.
Hence,
\begin{eqnarray*}
\sum_{j\in K} g_ju_j = \sum_{x\in
X} n_xw_{x1}\ldots w_{xz_x},
\end{eqnarray*}
where $n_x\in F_K\cap N,\,w_{pq}\in F_K,\,z_x\geqslant 0$.\\
We define element $v$ by $v = \sum_{x\in X} [\ldots[n_xw_{x1}]\ldots
w_{xz_x}]$. Then $v\in F_K\cap N$ and
\begin{eqnarray}\label{alg1_lm1_2}
D_j(v) \equiv \sum_{x\in X} D_j(n_x)w_{x1}\ldots
w_{xz_x}\mod{N_U}\,(j\in K).
\end{eqnarray}
We have
\begin{gather}
\sum_{j\in K} g_j\sum_{x\in X} D_j(n_x)w_{x1}\ldots w_{xz_x} =\notag\\
\sum_{x\in X}\sum_{j\in K}
g_jD_j(n_x)w_{x1}\ldots w_{xz_x} =\notag\\
\sum_{x\in X}n_xw_{x1}\ldots w_{xz_x}=\sum_{j\in K} g_ju_j.\notag
\end{gather}
Thus,
\begin{eqnarray}\label{alg1_lm1_3}
\sum_{x\in X} D_j(n_x)w_{x1}\ldots w_{xz_x}=u_j\,(j\in K).
\end{eqnarray}
From (\ref{alg1_lm1_2}), (\ref{alg1_lm1_3}) follows $D_j(v) \equiv
u_j\mod{N_U}\,(j\in K)$.
\end{proof}

\begin{lemma}\label{alg1_tm1}
Suppose $F$ is a free Lie algebra on free generators $\{g_j \mid  j\in J\}$, $N$ --- an ideal of $F$,
$D_j\,(j\in J)$ --- the Fox derivatives of $U(F)$, $K\subseteq J$, $F_K$
the subalgebra of $F$, generated by $\{g_j\mid  j\in K\}$, $\{u_j\mid j\in K\}$ an
elements of $U(F)$ such that $u_j$ is equal to zero for all but
a finite number of $j$. If
\begin{eqnarray}\label{alg1_tm1_1}
\sum_{j\in K} g_ju_j\equiv 0\mod{N_U},
\end{eqnarray}
then there is an element $v\in \mbox{\rm id}_F (F_K\cap N)$  such that $D_j(v)\equiv u_j \mod{N_U}$ $(j\in K)$.
\end{lemma}
\begin{proof}
Let $U(F)$ has for basis the set of monomials of the form
(\ref{alg1_1}).
There exist a monomials $f_1,\ldots, f_m$
$(f_i\neq f_j$ if $i\neq j)$ of the form (\ref{alg1_1}) with $\mu=\nu=\eta=0$ such that
\begin{eqnarray}\label{alg1_tm1_2}
u_j \equiv \sum_{l=1}^m u_{jl}f_l\mod{N_U}\,(j\in K),
\end{eqnarray}
where $u_{jl}\in U(F_K)$.
From (\ref{alg1_tm1_1}) and
\begin{eqnarray*}
\sum_{j\in K} g_ju_j \equiv
\sum_{l=1}^m  \sum_{j\in K}
g_ju_{jl}f_l\mod{N_U},
\end{eqnarray*}
follows
\begin{eqnarray}\label{alg1_tm1_3}
\sum_{j\in K} g_ju_{jl} \equiv
0 \mod{N_U}\,(1\leqslant l\leqslant m).
\end{eqnarray}
From (\ref{alg1_tm1_3}) and lemma \ref{alg1_lm1}
follows the existence of an elements $v_1,\ldots, v_m$ of $F_K\cap N$ such that
\begin{eqnarray}\label{alg1_tm1_4}
D_j(v_l) \equiv u_{jl}\mod{N_U}\,(j\in K).
\end{eqnarray}
Let $f_l=d_{l1}\ldots d_{lz_l}$. We define element $v$ by $v =
\sum_{l=1}^m [\ldots[v_ld_{l1}]\ldots d_{lz_l}]$.\\
Then
$v\in\mbox{\rm id}_F (F_K\cap N)$ and
\begin{eqnarray}\label{alg1_tm1_5}
D_j(v) \equiv \sum_{l=1}^m D_j(v_l)f_l\mod{N_U}\,(j\in K).
\end{eqnarray}
From (\ref{alg1_tm1_4}), (\ref{alg1_tm1_5}) follows
\begin{eqnarray}\label{alg1_tm1_6}
D_j(v) \equiv \sum_{l=1}^m u_{jl}f_l\mod{N_U}\,(j\in K).
\end{eqnarray}
Using (\ref{alg1_tm1_2}), (\ref{alg1_tm1_6}) we now find $D_j(v)
\equiv u_j\mod{N_U}$ $(j\in K)$.
\end{proof}

\begin{theorem}\label{alg1_tm2}
Suppose $F$ is a free Lie algebra on free generators $\{g_j \mid  j\in J\}$, $N$ --- an ideal of $F$,
$D_j~(j\in J)$ --- the Fox derivatives of $U(F)$, $K\subseteq J$, $F_K$
 --- the subalgebra of $F$, generated by $\{g_j\mid  j\in K\}$, $v\in F$.
Then
\begin{eqnarray}\label{alg1_tm2_1}
D_k(v)\equiv 0 \mod{N_U} \,(k \in J\verb|\|K)
\end{eqnarray}
if and only if there are an elements $v_0\in F_K$ and $v_1\in \mbox{\rm id}_F (F_K\cap N)$ such that $v\equiv v_0 + v_1
\mod{[N,N]}$.
\end{theorem}
\begin{proof}{\bf Sufficiency} is obvious.

{\bf Necessity}. Let $U(F)$ has for basis the set of monomials of
the form (\ref{alg1_1}). From (\ref{alg1_tm2_1}) follows that if $j\in J\verb|\|K$, then $g_jD_j(u)$ --- a linear combination of monomials of the form (\ref{alg1_1}) with $\mu+\eta \geqslant 1$. From $j\in K$ follows that $g_jD_j(u)$ --- a linear combination of monomials of the form (\ref{alg1_1}) with $\mu+\nu \geqslant 1$. Thus $v$ --- a linear combination of monomials of the form (\ref{alg1_1}) with $\mu+\nu+\eta = 1$, whence $v\in F_K+N$.
Then there is an element $v_0\in F_K$ such that $v-v_0\in N$. From $D_k(v-v_0)\equiv 0 \mod{N_U}$ $(k\in J\verb|\|K)$ and $v-v_0=\sum_{j\in J} g_jD_j(v-v_0)$ follows
\begin{eqnarray}\label{alg1_tm2_2}
\sum_{j\in K}
g_jD_j(v-v_0)\equiv 0 \mod{N_U}.
\end{eqnarray}
From (\ref{alg1_tm2_2}) and lemma \ref{alg1_tm1} follows that
there exists $v_1\in\mbox{\rm id}_F (F_K\cap N)$ such that
\begin{eqnarray*}
D_k(v-v_0)\equiv D_k(v_1)\mod{N_U},\,k\in J.
\end{eqnarray*}
Then
\begin{eqnarray*}
D_k(v-v_0 - v_1)\equiv \,0\mod{N_U},\,k\in J,
\end{eqnarray*}
whence $v-v_0 - v_1\in [N,N]$ \cite{Hm}.
\end{proof}

\section{A theorem on freedom for relatively free Lie algebras}
Let $F$ be a free Lie algebra,  $H$ --- a subalgebra of $F$; $N=N_1 \geqslant \ldots \geqslant N_t \geqslant \ldots$ --- a descending series of ideals of $F$ such that $[N_i,N_j\,]\leqslant N_{i+j}$, $H_i=H\cap N_i$.

Let $\{e_i,\,i\in I\}$ be an ordered basis in $N_l$,
$\{e_i,\,i\in I\}\cup\{a_{l-1,i},\,i\in I_{l-1}\}$ --- an ordered basis in
$N_l+H_{l-1}\,(\{a_{l-1,i},\,i\in I_{l-1}\}\subseteq H_{l-1})$, $\{e_i,\,i\in I\}\cup\{a_{l-1,i},\,i\in I_{l-1}\}\cup\{b_{l-1,i},\,i\in J_{l-1}\}$ --- an ordered basis in $N_{l-1}$ and iterate this procedure we obtain
an ordered basis $M=\{e_i,\,i\in I\}\cup\{a_{l-1,i},\,i\in I_{l-1}\}\cup\ldots\cup\{b_{1,i},\,i\in J_1\}$ in $N$.

Let $M\cup\{c_i,\,i\in J\}\,(\{c_i,\,i\in J\}\subseteq H)$ be an ordered basis in $H+N$, $M\cup\{c_i,\,i\in J\}\cup\{d_i,\,i\in I_F\}$ --- an ordered basis in $F$.

We assume that $d_t<c_i<b_{k,j}<a_{s,m}<e_p$ and introduce on ${\bf N}\times {\bf N}$ the lexicographic ordering $(s,m)<(k,j)$ if $s<k$ or $s=k$ and $m<j$.
Then $U(F)$ has for basis the set $U_F$ of monomials of the form
\begin{eqnarray}\label{alg2_1}
d_{t_1}\ldots
d_{t_\theta}c_{i_1}\ldots
c_{i_\eta}b_{j_1,k_1}\ldots b_{j_\nu,k_\nu}a_{s_1,m_1}\ldots a_{s_\mu,m_\mu}e_{p_1}\ldots e_{p_\gamma},
\end{eqnarray}
where $t_1\leqslant\ldots\leqslant t_\theta,\,i_1\leqslant\ldots\leqslant i_\eta,\,(j_1,k_1)\leqslant\ldots\leqslant (j_\nu,k_\nu)$, $(s_1,m_1)\leqslant\ldots\leqslant (s_\mu,m_\mu),\,p_1\leqslant\ldots\leqslant p_\gamma$;
$\theta,\,\eta,\,\nu,\,\mu,\,\gamma\geqslant 0$.

We now introduce an ordering on the set $U_F$ 
by taking words of greater
length to be greater than words of smaller length, and using
lexicographic order for words of equal length (from left to right).

If $u$ --- monomials of the form (\ref{alg2_1}) then we denote $\theta+\eta+\nu+\mu+\gamma$ by $d(u)$.
Let $u,\,v\in U_F$, $uv=\sum_{i={1}}^{p}\alpha_i m_i$ where $0\neq\alpha_i\in P,\,i=1,\ldots,p$, $m_i\in U_F,\,i=1,\ldots,p$
and $m_i<m_j$ for $i<j$. It is clear that $d(m_p)= d(u)+d(v)$.
From (\ref{sob_pr}) we obtain $d(m_i)\leqslant d(u)+d(v)-1,\,i=1,\ldots,p-1$.
We denote $m_p$ by $u\circ v$. It can be verified directly that if $u_1$, $u_2$, $v_1$, $v_2$ --- an elements of $U_F$ and $u_1<u_2$, $v_1\leqslant v_2$ then
$u_1\circ v_1<u_2\circ v_2$.

We denote the set of monomials of the form (\ref{alg2_1}) with  $\theta> 0$ and $\nu = \mu = \gamma = 0$ by $S_\beta$,
the set of monomials of the form (\ref{alg2_1}) with $\theta = \nu =\mu = \gamma =0$ by $S_\alpha$,
$S_\alpha\cup S_\beta$ by $S$; the ideal in $U(N)$ generated by all products of the form $N_{i_1}\cdots N_{i_t}$ $(i_1 + \cdots + i_t\geqslant i)$ by $\Delta_i$; the ideal in $U(H_1)$ generated by all products of the form $H_{i_1}\cdots H_{i_t}$ $(i_1 + \cdots + i_t\geqslant i)$ by $\Delta_i^\prime$.

By $\Delta_0$, $\Delta_0^\prime$ we denote $U(N)$, $U(H_1)$ respectively.
It is clear that $\Delta_1=U_0(N)$, $\Delta_1^\prime=U_0(H_1)$.
Since $[N_i,N_j\,]\leqslant N_{i+j}$ it follows that any element of $S\Delta_i\setminus (S\Delta_{i+1}+(N_l)_U)$ is
a linear combination of linearly independent over $P$ modulo $S\Delta_{i+1}+(N_l)_U$ monomials of the form (\ref{alg2_1})
with $j_1+\ldots + j_\nu+s_1+\ldots + s_\mu= i$ and $\gamma=0$.

Let $u$ --- a linear combination of linearly independent over $P$ modulo $S\Delta_{i+1}+(N_l)_U$ monomials $u_1<\ldots<u_n$ of the form (\ref{alg2_1}) with $j_1+\ldots + j_\nu+s_1+\ldots + s_\mu= i$ and $\gamma=0$;
$v$ --- a linear combination of linearly independent over $P$ modulo $S\Delta_{j+1}+(N_l)_U$ monomials $v_1<\ldots<v_m$ of the form (\ref{alg2_1}) with $j_1+\ldots + j_\nu+s_1+\ldots + s_\mu= j$ and $\gamma=0$. Then
$u_n\circ v_m$ --- monomial with $j_1+\ldots + j_\nu+s_1+\ldots + s_\mu= i+j$ and $\gamma=0$. Hence if $u\in S\Delta_i\setminus S\Delta_{i+1}\mod{(N_l)_U}$, $v\in S\Delta_j\setminus S\Delta_{j+1}\mod{(N_l)_U}$,
then $uv\in S\Delta_{i+j}\setminus S\Delta_{i+j+1}\mod{(N_l)_U}$.

We now show that
\begin{eqnarray}\label{fm2_pr_gr}
U(H_1)\cap \Delta_i \equiv \Delta_i^\prime\mod{(N_l)_U}.
\end{eqnarray}
It is clear that $\Delta_i^\prime\subseteq U(H_1)\cap \Delta_i$.
Since any element of $U(H_1)\cap \Delta_i$ is a linear combination of monomials
of the form (\ref{alg2_1}) with $\gamma > 0$ or $\theta = \eta = \nu = \gamma=0$ and $s_1+\ldots + s_\mu\geqslant i$ it follows that this element belongs to $\Delta_i^\prime+(N_l)_U$.

\begin{lemma}\label{lm1_2}
Suppose $F$ is a free Lie algebra, 
$H$ --- a subalgebra of $F$; 
$N=N_1 \geqslant \ldots \geqslant N_z \geqslant \ldots $ --- a descending series of ideals of $F$ such that $[N_i,N_j\,]\leqslant N_{i+j}$; $U(F)$ has for basis the set of monomials of the form
(\ref{alg2_1}); $f_1,\ldots,f_l\, (f_i\neq f_j$ for $i\neq j)$, $g_1,\ldots,g_k\, (g_i\neq g_j$ for $i\neq j)$ --- an elements of $S$; $f_ig_j=W_{ij}$, where $W_{ij}$ --- a linear combination of linearly independent over $P$ monomials of the form (\ref{alg2_1});
$\{f_1,\ldots,f_l,g_1,\ldots,g_k\}\not\subseteq S_\alpha$. Then there exist
$M\in S_\beta,\,i_0,\,j_0$ such that $M$ --- monomial of $W_{i_0j_0}$; monomials of $W_{ij}$ are all different from $M$ for $(i_0,j_0)\neq (i,j)$.
\end{lemma}
\begin{proof}
If $u$ --- monomial of the form (\ref{alg2_1}) then we set $d_\theta(u)=\theta$.

We denote $\max\,(d_\theta(f_1),\ldots,d_\theta(f_l))$ by $x$; $\max\,(d_\theta(g_1),\ldots,d_\theta(g_k))$ by $y$.
By the assumptions of the lemma $x+y>0$.

Set $f_{i_0}=\max\,\{f_i\mid d_\theta(f_i)=x\}$; $g_{j_0}=\max\,\{g_j\mid d_\theta(g_j)=y\}$; $M=f_{i_0}\circ g_{j_0}$.
Since $d_\theta(M)=d_\theta(f_{i_0})+d_\theta(g_{j_0})=x+y$ we have $M\in S_\beta$. It is clear that for $(i_0,j_0)\neq (i,j)$ monomials of $W_{ij}$ are all different from $M$.
\end{proof}

\begin{proposition}\label{tm4}
Suppose $F$ is a free Lie algebra on free generators $y_1,\ldots,y_n$, $n\geqslant 3$, $N_{11}$  ---  an ideal of $F$, $F/N_{11}$  ---  a relatively free Lie algebra,
\begin{eqnarray}\label{tm4_0}
N_{11} \geqslant \ldots \geqslant N_{1,m_1+1}=N_{21} \geqslant \ldots \geqslant N_{s,m_s+1},
\end{eqnarray}
where $N_{kl}$ --- the $l$-th term of the lower central series of $N_{k1}$.\\
Let $R$ --- be an ideal of $F$, $R\leqslant N_{11}$; $H$ --- the subalgebra of $F$, generated by $y_1,\ldots,y_{n-1}$; $j\in \{1,\ldots,m_1+1\}$ such that $H\cap (R+N_{1j})\neq H\cap N_{1j}$.\\
If $(k,l)\geqslant (1,j)$ then $H\cap (R+N_{kl})\neq H\cap N_{kl}$.
\end{proposition}
\begin{proof}
It is clear that $H\cap N_{mt}\supseteq (H\cap N_{m1})_{(t)}$.
Let $\phi$ be the endomorphism of $F$ defined by
$\phi(y_n)=0$, $\phi(y_i)=y_i$, $i=1,\ldots,n-1$.
Then $\phi(N_{mt})=(H\cap N_{m1})_{(t)}$, therefore  $u=\phi(u)\in (H\cap N_{m1})_{(t)}$ for any $u\in H\cap N_{mt}$.
Thus $H\cap N_{mt}=(H\cap N_{m1})_{(t)}$.

We denote by $N$ the algebra $H\cap N_{11}$.
Let $\{x_z \mid z\in I\}\,(I\subseteq {\bf N})$ be a free set of generators of $N$,
$\{\partial_z \mid z\in I\}$ --- the Fox derivatives of $N$.
Since free set of generators of $H$ contains more than one element and
$N$ --- an ideal of $H$ it follows that free set of generators of $N$ contains more than one element.
So we may assume that $|I| > 1$.

We now show that
\begin{eqnarray}\label{pr1}
H\cap (R+N_{1i})> N_i,\,i=j,\ldots,m_1+1.
\end{eqnarray}
If $j=m_1+1$ the truth of the statement is obvious.
Assume inductively that $H\cap (R+N_{1i})> N_i\,(i=j,\ldots, l;\,l\leqslant m_1)$.
We need to show that $H\cap (R+N_{1,l+1})> N_{l+1}$.

It is well known that $U(N)$ --- a free associative algebra with a free set $\{x_z  \mid z\in I\}$ of generators.
Hence $v\in N_{(l)}$ if and only if $v\in U_0(N)^l$ for any $v\in N$.
Since $v=\sum_{z\in I }x_z\partial_z(v)$ it follows that
$v\in N_{(l)}$  if and only if $\partial_z(v)\in U_0(N)^{l-1},\,z\in I$.

Let $v\in (H\cap (R+N_{1l}))\setminus N_{(l)}$. Then $v\notin U_0(N)^l$; further we may assume that
$\partial_1(v)\notin U_0(N)^{l-1}$ without loss of generality.
We denote by $w$ the element $[v,x_2]$.

From $v\in H\cap (R+N_{1l})$ follows that $w\in H\cap (R+N_{1,l+1})$.
Since $\partial_1(w)=\partial_1(v)x_2$ we have $\partial_1(w)\notin U_0(N)^l$. Therefore $w\notin U_0(N)^{l+1}$,
whence $w\in (H\cap (R+N_{1,l+1}))\setminus N_{(l+1)}$ and, by induction on $l$,
$H\cap (R+N_{1l})> N_l$ for each term $N_{1l}\,(l\geqslant j)$ of series {\rm (\ref{tm4_0})}.
From {\rm (\ref{pr1})} follows that $H\cap (R+N_{21})> H\cap N_{21}$.

Since $H\cap (R+N_{kl})\geqslant (H\cap (R+N_{k1}))_{(l)}$ and $H\cap N_{kl}=(H\cap N_{k1})_{(l)}$
we reach the conclusion by noting that if $H\cap (R+N_{k1})> H\cap N_{k1}$ then
$(H\cap (R+N_{k1}))_{(l)}> (H\cap N_{k1})_{(l)}$.
\end{proof}

\begin{proposition}\label{tm2}
Suppose $F$ is a free Lie algebra on free generators $y_1,\ldots,y_n$,
$H$ --- the subalgebra of $F$, generated by $y_1,\ldots,y_{n-1}$, $N_{11}$  --- an ideal of $F$, $F/N_{11}$  ---  a relatively free Lie algebra,
\begin{eqnarray}\label{tm2_0}
N_{11} \geqslant \ldots \geqslant N_{1,m_1+1}=N_{21} \geqslant \ldots \geqslant N_{s,m_s+1},
\end{eqnarray}
where $N_{kl}$ --- the $l$-th term of the lower central series of $N_{k1}$.\\
Let $r\in N_{1i}\backslash N_{1,i+1}\,(i\leqslant m_1)$, $R = \mbox{\rm id}_F (r)$.
If $H\cap (R+N_{21})=H\cap N_{21}$ then\\
$H\cap (R+N_{kl})=H\cap N_{kl}\,(k> 1)$ for each term $N_{kl}\,(k> 1)$ of series {\rm (\ref{tm2_0})}.
\end{proposition}
\begin{proof}
If $n=2$ then $H\cap N_{21}=0$. By the assumptions of the proposition $H\cap N_{21}=H\cap (R+N_{21})$ hence
$H\cap (R+N_{21})=0$ and $H\cap (R+N_{kl})=H\cap N_{kl}=0\,(k> 1)$.
So we may assume that $n>2$.

We denote by $D_1,\ldots,D_n$ the Fox derivatives of $F$.
Let us now suppose that $D_n(r)\equiv 0\mod {(R+N_{21})_U}$.
Theorem~\ref{alg1_tm2} tells us that there exist an elements $r_0\in H$ and $r_1\in \mbox{\rm id}_F (H\cap (R+N_{21}))$ such that $r\equiv r_0 + r_1 \mod{(R+N_{21})_{(2)}}$.
Since $r$, $r_1$ --- an elements of $R+N_{21}$ we have $r_0\in H\cap (R+N_{21})$, $r\in \mbox{\rm id}_F (H\cap (R+N_{21}))+(R+N_{21})_{(2)}$.
From $H\cap (R+N_{21})=H\cap N_{21}$ follows $r\in N_{21}+(R+N_{21})_{(2)}$. Hence $r\in N_{1,i+1}$, a contradiction.
So we may assume that $D_n(r)\not\equiv 0\mod {(R+N_{21})_U}$.

Now assume inductively that $H\cap (R+N_{ij}) = H\cap N_{ij}\,(i=2,\ldots, k;\,j=1,\ldots, l;\,l\leqslant m_k)$.
We need to show that $H\cap (R+N_{k,l+1})=H\cap N_{k,l+1}$.

Let us denote the ideal $R+N_{k1}$ by $N$ and denote by $N_m$ the ideal $R+N_{km}$, $m\in {\bf N}$.
It is not hard to verify that $[N_p\,,N_q\,]\leqslant N_{p+q}$.
Then $U(F)$ has for basis the set of monomials of the form
(\ref{alg2_1}), $S=S_\alpha\cup S_\beta$.

Let $\{x_{kz} \mid  z\in {\bf N}\}$ be a free set of generators of $N$;
$\{\partial_{kz} \mid  z\in {\bf N}\}$ ---  the Fox derivatives of $U(N)$;
$v\in H\cap N_{l+1}$.
From $N_{l+1}=R+N_{k,l+1}=R+N_{(l+1)}$ follows the existence of an element $u\in N_{(l+1)}$ such
that $v-u\in \mbox{\rm id}_F (r)$.
From (\ref{four}) follows the existence of an elements $k_p\in U(N)$, $f_p\in S$, $p=1,\ldots,d$ $(f_p\neq f_j\mbox{ for }p\neq j)$ such that $D_m(v-u)\equiv D_m(r)\cdot  \sum_{p={1}}^{d} f_p k_p\mod{(N_l)_U},\,m=1,\ldots,n$.

By (\ref{alg1_dcf}), $D_m(u)=\sum_{z\in {\bf N}} D_m(x_{kz})\partial_{kz}(u)$ thus
\begin{eqnarray}\label{tm2_3}
D_m(v)\equiv D_m(r)\cdot  \sum_{p={1}}^{d} f_p k_p +\sum_{z\in {\bf N}} D_m(x_{kz})\partial_{kz}(u)\mod{(N_l)_U},
\end{eqnarray}
$m=1,\ldots,n$. From $v\in H$ and (\ref{tm2_3}) follows
\begin{eqnarray}\label{tm2_4}
0\equiv D_n(r)\cdot  \sum_{p={1}}^{d} f_p k_p+ \sum_{z\in {\bf N}} D_n(x_{kz})\partial_{kz}(u)\mod{(N_l)_U}.
\end{eqnarray}
From $D_n(r)\not\equiv 0\mod {N_U}$ follows the existence of an elements
$0\neq\gamma_p\in P$, $g_p\in S$, $p=1,\ldots,q;\,g_p\neq g_j$ for $p\neq j$, such that $D_n(r) \equiv  \sum_{p={1}}^q  \gamma_p g_p \mod {N_U}$.

Let us now suppose that there exist $l_0,\,j$
such that $l_0< l$; $k_j\in U_0(N)^{l_0}\setminus U_0(N)^{l_0+1}\mod{(N_l)_U}$; $k_p\in U_0(N)^{l_0}\mod{(N_l)_U}$, $p=1,\ldots,d$.
We denote $\max\,(g_1,\ldots,g_q)$ by $g_a$; $\max\,\{f_i\mid k_i\in U_0(N)^{l_0}\setminus U_0(N)^{l_0+1}\mod{(N_l)_U}\}$ by $f_b$;
$g_a\circ f_b$ by $M$. We have
\begin{eqnarray}\label{tm2_4_1}
D_n(r)\cdot  \sum_{p={1}}^{d} f_p k_p \equiv \gamma_a M k_b + h  \mod{(N_l)_U},
\end{eqnarray}
where $h$ --- a linear combination of monomials of the form (\ref{alg2_1}) and if this monomials are not in $S\cdot U_0(N)^{l_0 +1}$ then this monomials are not in $M U(N)$.

Since $u\in N_{(l+1)}$ it follows that $\sum_{z\in {\bf N}} D_m(x_{kz})\partial_{kz}(u)\in S\cdot U_0(N)^l$,
therefore (\ref{tm2_4_1}) contradicts (\ref{tm2_4}). So we may assume that
\begin{eqnarray}\label{tm2_4_1_1}
k_p\in U_0(N)^l\mod{(N_l)_U},\,p=1,\ldots,d.
\end{eqnarray}
By (\ref{tm2_3}), (\ref{tm2_4_1_1}) we obtain
\begin{eqnarray}\label{tm2_6}
D_m(v)\equiv \sum_{i={1}}^{d_m} f_{im}v_{im}  \mod{(N_l)_U},~m=1,\ldots,n,
\end{eqnarray}
where $f_{im}\in S$, $v_{im}\in U_0(N)^l$.

Let $\{x_z \mid  z\in I\subseteq {\bf N})\}$ be a free set of generators of $H\cap N$;
$\{\partial_z \mid  z\in I\}$ ---  the Fox derivatives of $U(H\cap N)$.
If $D_j(x_z) \equiv 0\mod{N_U}\,(j=1,\ldots,n-1)$ then $D_j(x_z) \equiv 0\mod{U(H) (H\cap N)},\,j=1,\ldots,n-1$;
hence $x_z\in [H\cap N,H\cap N]$ \cite{Hm}.
Since $x_z\notin  [H\cap N,H\cap N]$ it follows that there exist $j_z\in \{1,\ldots,n-1\}$ such that $D_{j_z}(x_z) \not\equiv 0\mod{N_U}$.
Thus $D_{j_z}(x_z) \equiv \gamma_z M_z+V_z\mod{N_U}$,
where $0\neq\gamma_z\in P$; $M_z\in S_\alpha$; $V_z$ --- a linear combination of linearly independent over $P$ modulo $N_U$ monomials of the form (\ref{alg2_1}); this monomials belong to $S_\alpha$ and all different from $M_z$.
It is not hard to verify that if $t\neq z$ then there exists $\gamma\in P$ such that putting $x_t^\prime=x_t-\gamma x_z$ we obtain $D_{j_z}(x_t^\prime) \equiv W\mod{N_U}$ where $W$--- a linear combination of linearly independent over $P$ modulo $N_U$ monomials of the form (\ref{alg2_1}); this monomials belong to $S_\alpha$ and all different from $M_z$.

It is no restriction to assume that $v\in (x_1,\ldots,x_a)$;
for any $m\in \{1,\ldots,a\}$ there exists $j_m\in \{1,\ldots,n\}$ such that
$D_{j_m}(x_m) \equiv \gamma_m\cdot M_m+V_m \mod{N_U}$ and if $z\neq m$ then $M_z$ and monomials of $V_z$ are all different from $M_m$.

From $v\in H\cap N_l=H\cap N_{kl}$ follows $v\in (H\cap N_{k1})_{(l)}\subseteq (H\cap N)_{(l)}$. Hence
\begin{eqnarray}\label{tm2_6_0_6}
\partial_z(v)\in U_0(N)^{l-1},\,z=1,\ldots,a.
\end{eqnarray}
Let us now show that $\partial_z(v)\in U_0(N)^l\mod{(N_l)_U},\,z=1,\ldots,a$.
Suppose that there exists $i\in \{1,\ldots,a\}$ such that $\partial_i(v)\not\in U_0(N)^l\mod{(N_l)_U}$.

By $D_{j_i}(v)= \sum_{z=1}^a D_{j_i}(x_z)\partial_z(v)$ and (\ref{tm2_6_0_6}) we obtain
\begin{eqnarray}\label{tm2_6_0}
D_{j_i}(v) \equiv \gamma_i\cdot M_i\partial_i(v)+W \mod{(N_l)_U}.
\end{eqnarray}
where $W$ --- a linear combination of monomials of the form (\ref{alg2_1}) and if this monomials are not in $S\cdot U_0(N)^l$ then this monomials are not in $M_i U(N)$.

Thus (\ref{tm2_6_0}) contradicts (\ref{tm2_6}) hence
\begin{eqnarray}\label{tm2_6_00}
\partial_z(v)\in  U_0(N)^l\mod{(N_l)_U} \,(z=1,\ldots,a)
\end{eqnarray}
as we wished to show.

Since $U(H_1)\cap \Delta_l\equiv \Delta_l^\prime\mod{(N_l)_U}$ it follows from (\ref{tm2_6_00}) that $\partial_z(v)\in \Delta_l^\prime\mod{(N_l)_U}$. Therefore there exists $v_z\in \Delta_l^\prime$  such that
$\partial_z(v)-v_z\in U(H_1)\cap (N_l)_U$ whence $\partial_z(v)\in U(H_1)(H\cap N_l) + \Delta_l^\prime$ $(z=1,\ldots,a)$.

We have, by the inductive assumption, that $H\cap N_l)=H\cap N_{kl}=(H\cap N_{k1})_{(l)}\subseteq U_0(H_1)^l$, whence $\partial_z(v)\in \Delta_l^\prime$ $(z=1,\ldots,a)$.

If $t\leqslant l$ then $H_t=H\cap N_{kt} = (H_1)_{(t)}\subseteq U_0(H_1)^t$ hence
$\Delta_l^\prime\subseteq U_0(H_1)^l$, i.e. $\partial_z(v)\in U_0(H_1)^l$, $z=1,\ldots,a$.
Thus $v\in U_0(H_1)^{l+1}$ whence $v\in (H\cap N_{k1})_{(l+1)}=H\cap N_{k,l+1}$ and, by induction on $l$,
$H\cap (R+N_{kl})=H\cap N_{kl}$ for each term $N_{kl}\,(k> 1)$ of series {\rm (\ref{tm2_0})}.
\end{proof}

\begin{proposition}\label{tm5}
Suppose $F$ is a free Lie algebra on free generators $y_1,\ldots,y_n$,
$H$ --- the subalgebra of $F$, generated by $y_1,\ldots,y_{n-1}$, $N$ --- an ideal of $F$, $F/N$  ---  a relatively free Lie algebra,
\begin{eqnarray}\label{tm5_0}
N=N_{11} \geqslant \ldots \geqslant N_{1,m_1+1}=N_{21} \geqslant \ldots \geqslant N_{s,m_s+1},
\end{eqnarray}
where $N_{kl}$ --- the $l$-th term of the lower central series of $N_{k1}$.
Let $r\in N_{1i}\backslash N_{1,i+1}\,(i\leqslant m_1)$, $R = \mbox{\rm id}_F (r)$.
Then $H\cap (R+N_{1l})=H\cap N_{1l}$ for each term $N_{1l}$ of series {\rm (\ref{tm5_0})}.
\end{proposition}
\begin{proof}
We denote by $N_m$ the ideal $R+N_{1m}$, $m\in {\bf N}$. It is not hard to verify that $[N_p\,,N_q\,]\leqslant N_{p+q}$.
Then $U(F)$ has for basis the set of monomials of the form
(\ref{alg2_1}), $S=S_\alpha\cup S_\beta$.

Let $\{x_z \mid  z\in I\subseteq {\bf N})\}$ be a free set of generators of $N$;
$\{\partial_z \mid  z\in I\}$ ---  the Fox derivatives of $U(N)$.
Free generator $x_z$ we name $\alpha$-generator if $x_z\in H\cap N\mod{[N,N]}$.
It is clear that if $x_z$ --- $\alpha$-generator then $D_n (x_z)\equiv \,0\mod{N_U}$ and
$D_j (x_z)\in U(H)\mod{N_U},\,j\in \{1,\ldots,n-1\}$.

Let us now suppose that $x\in \{x_z \mid  z\in I\subseteq {\bf N})\}$; $D_n (x)\equiv \,0\mod{N_U}$ and
$D_j (x)\in U(H)\mod{N_U},\,j\in \{1,\ldots,n-1\}$.
Since $x=\sum_{j=1}^{n-1}y_jD_j(x)\equiv 0 \mod{N_U}$,
by Lemma~\ref{alg1_lm1} we obtain the existence of an element $y\in H\cap N$ such that
$D_j(x)\equiv D_j (y) \mod{N_U},\,j\in \{1,\ldots,n-1\}$.
Therefore, $x-y\in [N,N]$ \cite{Hm}, i.e. $x$ --- $\alpha$-generator.

Thus if $x_z$ is not $\alpha$-generator then $D_n (x_z)\not\equiv \,0\mod{N_U}$ or $D_n (x_z)\equiv \,0\mod{N_U}$ and
there exists $j\in \{1,\ldots,n-1\}$ such that $D_j (x_z)\notin U(H)\mod{N_U}$.

Since $x_z\notin  [N,N]$ it follows that there exist $j_z\in \{1,\ldots,n\}$ such that $D_{j_z}(x_z) \not\equiv 0\mod{N_U}$ \cite{Hm}.
Thus $D_{j_z}(x_z) \equiv \gamma_z M_z+V_z\mod{N_U}$,
where $0\neq\gamma_z\in P$; $M_z\in S$; $V_z$--- a linear combination of linearly independent over $P$ modulo $N_U$ monomials of the form (\ref{alg2_1}); this monomials belong to $S$ and all different from $M_z$. It is clear that if $t\neq z$ then there exists $\gamma\in P$ such that putting $x_t^\prime=x_t-\gamma x_z$ we obtain $D_{j_z}(x_t^\prime) \equiv W\mod{N_U}$ where $W$--- a linear combination of monomials of the form (\ref{alg2_1}); this monomials belong to $S$ and all different from $M_z$.

It is no restriction to assume that $r,\,v\in (x_1,\ldots,x_a)$;
for any $m\in \{1,\ldots,a\}$ there exists $j_m\in \{1,\ldots,n\}$ such that
$D_{j_m}(x_m) \equiv \gamma_m\cdot M_m+V_m \mod{N_U}$; if $x_m$ is not $\alpha$-generator and $D_n(x_m)\equiv \,0\mod{N_U}$ then $M_m\in S_\beta$; if $j>m$ then $M_j$ and monomials of $V_j$ are all different from $M_m$.

It is clear that $H\cap N_i=H\cap N_{1i}$.
Now assume inductively that $H\cap N_j=H\cap N_{1j}$ $(j=i,\ldots,l;\,l\leqslant m_1)$.
We need to show that $H\cap N_{l+1}=H\cap N_{1,{l+1}}$.
Let $v\in H\cap N_{l+1}$. From $N_{l+1}=R+N_{1,l+1}=R+N_{(l+1)}$ follows the existence of an element $u\in N_{(l+1)}$ such
that $v-u\in \mbox{\rm id}_F (r)$.
From (\ref{four}) follows the existence of an element $C\in U(F)$ such that $D_j(v-u)\equiv D_j(r)\cdot  C,\,j=1,\ldots,n$.

By (\ref{alg1_dcf}), $D_m(u)=\sum_{z\in I} D_m(x_z)\partial_z(u)$ thus
\begin{eqnarray}\label{tm3_3}
D_j(v)\equiv D_j(r)\cdot  C +\sum_{z\in I} D_j(x_z)\partial_z(u)\mod{(N_l)_U},
\end{eqnarray}
$j=1,\ldots,n$.
From $v\in H$ and (\ref{tm3_3}) follows
\begin{eqnarray}\label{tm3_4}
0\equiv D_n(r)\cdot  C + \sum_{z\in I} D_n(x_z)\partial_z(u)\mod{(N_l)_U}.
\end{eqnarray}

Since $r\in N_{(i)}$ it follows that $\partial_z(r)\in U_0(N)^{i-1},\,z\in I$.
Taking into account that $\Delta_j=U_0(N)^j\,(j\leqslant i)$ we may assert that from $r\notin  N_{(i+1)}$ follows the existence of $m$ such that $\partial_m(r)\in \Delta_{i-1}\setminus \Delta_i$ and $\partial_j(r)\in \Delta_i\, for\, j<m$.

By (\ref{alg1_dcf}), $D_{j_m}(r)= \sum_{z=1}^a D_{j_m}(x_z)\partial_z(r)$, hence
\begin{eqnarray}\label{tm2_6_01}
D_{j_m}(r) \equiv \gamma_m\cdot M_m\partial_m(r)+W_1 \mod{U(F)\Delta_i},
\end{eqnarray}
where $W_1$ --- a linear combination of monomials of the form (\ref{alg2_1}) and this monomials are not in $M_m U(N)$.

Since $v\in N_{(l)}$ it follows that $\partial_z(v)\in U_0(N)^{l-1},\,z\in I$.
Let us now show that $\partial_z(v)\in \Delta_l,\,z=1,\ldots,a$.
Suppose that there exists $k\in \{1,\ldots,a\}$ such that $\partial_k(v)\not\in \Delta_l$ and $\partial_j(v)\in \Delta_l\, for\, j<k$.

By (\ref{alg1_dcf}), $D_{j_k}(v)= \sum_{z=1}^a D_{j_k}(v)(x_z)\partial_z(r)$, hence
\begin{eqnarray}\label{tm2_6_01-1}
D_{j_k}(v) \equiv \gamma_k\cdot M_k\partial_k(v)+W_2\mod{U(F)\Delta_l},
\end{eqnarray}
where $W_2$ --- a linear combination of monomials of the form (\ref{alg2_1}) and this monomials are not in $M_k U(N)$.

Since $u\in N_{(l+1)}$ it follows that $\sum_{z\in I} D_m(x_z)\partial_z(u)\in SU_0(N)^l,\,z\in I$.
Then from $D_{j_k}(v)\notin U(F)\Delta_l$ and (\ref{tm3_3}) we obtain that
\begin{eqnarray}\label{tm3_4_1_1}
C \in U(F)\Delta_{l-i}\setminus U(F)\Delta_{l-i+1}.
\end{eqnarray}

We now show that if $x_z$ is not $\alpha$-generator then $\partial_z(r)\in U_0(N)^i$.
Let us suppose that there exists $b$ such that $x_b$ is not $\alpha$-generator; $\partial_b(r)\in U_0(N)^{i-1}\setminus U_0(N)^i$ and $\partial_c(v)\in U_0(N)^i$ for $x_c$ is not $\alpha$-generator and $c<b$.

Let $D_n(x_b)\not\equiv 0\mod{N_U}$.
If $x_c$ --- $\alpha$-generator, then $D_n(x_c)\equiv 0\mod{N_U}$, hence $D_n(x_c)\partial_c(r)\in U(F)\Delta_i$.
Then from $D_n(r)=\sum_z D_n(x_z)\partial_z(r)$ follows that
\begin{eqnarray}\label{tm3_4_2_alggr}
D_n(r)\equiv \gamma_b\cdot M_b\partial_b(r) +W_3 \mod{U(F)\Delta_i},
\end{eqnarray}
where $W_3$ --- a linear combination of monomials of the form (\ref{alg2_1});
this monomials are not in $M_b U(N)$. From (\ref{tm3_4_1_1}), (\ref{tm3_4_2_alggr}) we have $D_n(r)\cdot C\not\in U(F)\Delta_l$ which contradicts (\ref{tm3_4}).

Let $D_{j_b}(r)\equiv  \gamma_b\cdot M_b\partial_b(r) + W_4 \mod{U(F)\Delta_i}$;
$j_b\in \{1,\ldots,n-1\}$; $M_b\in S_\beta$; $W_4$ --- a linear combination of monomials of the form (\ref{alg2_1})
and this monomials are not in $M_b U(N)$.
By Lemma \ref{lm1_2}, taking into account that $D_m(v)\in U(H)\,(m=1,\ldots,n)$, we have $M_b\in S_\alpha$, a contradiction.

So we may assume that if $x_z$ is not $\alpha$-generator then $\partial_z(r)\in U_0(N)^i$.
We now show that if $x_z$ --- $\alpha$-generator then $\partial_z(r)\in U(N\cap H)\mod{U_0(N)^i}$.
Let us suppose that there exists $b$ such that $x_b$ --- $\alpha$-generator; $\partial_b(r)\not\in U(N\cap H)\mod{U_0(N)^i}$ and $\partial_c(r)\in U(N\cap H)\mod{U_0(N)^i}$ for $x_c$ --- $\alpha$-generator and $c<b$.
Then
\begin{eqnarray*}
D_{j_b}(r)\equiv  \gamma_b\cdot M_b\partial_b(r) + W_5 \mod{U(F)\Delta_i},
\end{eqnarray*}
where $M_b\in S_\alpha$, $W_5$ --- a linear combination of monomials of the form (\ref{alg2_1})
and this monomials are in $S_\alpha U(N)$; if some of this monomials are not in $S_\alpha U(N\cap H)$ then they are not in $M_b U(N)$.

If $m$ --- monomial of the form (\ref{alg2_1}) then we set $d_\nu(m)=j_1+\ldots + j_\nu$.
Set $d_\nu (\gamma_1m_1+\ldots+\gamma_xm_x)=\max\, (d_\nu(m_1),\ldots ,d_\nu(m_x))$,
where $\gamma_1,\ldots,\gamma_x$ --- non-zero elements of $P$, $m_1,\ldots,m_x$ --- monomials of the form (\ref{alg2_1}),
$m_i\neq m_j$ for $i\neq j$.

Monomial of the form (\ref{alg2_1}) with $\nu=\theta= 0$ we name $\alpha$-monomial.
From (\ref{sob_pr}) follows that if $u$, $v$ are $\alpha$-monomials then $uv$ --- a linear combination of $\alpha$-monomials.
Hence if $u$, $v$ are monomials of the form (\ref{alg2_1}) with $\theta= 0$ then
$d_\nu(uv)=d_\nu(u)+d_\nu(v)=d_\nu(u\circ v)$. It is not hard to verify that
if $u$, $v$ --- a linear combinations of monomials of the form (\ref{alg2_1}) with $\theta= 0$ then
$d_\nu(uv)=d_\nu(u)+d_\nu(v)$.

We denote by $\overline{D_{j_b}(v)}$ the linear combination of linearly independent over $P$ monomials of the form (\ref{alg2_1}) with $s_1+\ldots + s_\mu= l-1$ and $\theta=\nu= \gamma=0$ such that
 $D_{j_b}(v)\equiv \overline{D_{j_b}(v)}\mod{U(F)\Delta_l}$;
by $\overline{D_{j_b}(r)}$ the linear combination of linearly independent over $P$ monomials of the form (\ref{alg2_1}) with $j_1+\ldots + j_\nu+s_1+\ldots + s_\mu= i-1$ and $\theta=\gamma=0$ such that
 $D_{j_b}(r)\equiv \overline{D_{j_b}(r)}\mod{U(F)\Delta_i}$;
 by $\overline{C}$ the linear combination of linearly independent over $P$ monomials of the form (\ref{alg2_1}) with $j_1+\ldots + j_\nu+s_1+\ldots + s_\mu= l-i$ and $\theta=\gamma=0$ such that
 $C\equiv \overline{C}\mod{U(F)\Delta_{l-i+1}}$.

Thus we have 0=$d_\nu(\overline{D_{j_b}(v)})$ and $d_\nu(\overline{D_{j_b}(v)})=d_\nu( \overline{D_{j_b}(r)}\cdot  \overline{C})$
which contradicts $d_\nu( \overline{D_{j_b}(r)}\cdot  \overline{C})>0$.
This proves that if $k\in \{1,\ldots,a\}$ then $\partial_k(r)\in U_0(N)^i$ or $\partial_k(r)\in U(H\cap N)\mod{U_0(N)^i}$ and $x_k\in H\cap N\mod{[N,N]}$.
Let $\psi$ be the endomorphism of $U(F)$ defined by
$\psi(y_n)=0$, $\psi(y_i)=y_i$, $i=1,\ldots,n-1$.
We denote by $\hat{r}$ the element $\psi(r)$. It is clear that $\hat{r}\in H\cap N$.

Since $r= \sum_{t={1}}^a x_t\partial_t(r)$, we have
\begin{eqnarray}\label{main_009-1}
\hat{r}= \sum_{t={1}}^a \psi(x_t)\psi(\partial_t(r)),
\end{eqnarray}
where $\partial_t(r)$, $\psi(\partial_t(r))\in U_0(N)^i$ or
\begin{eqnarray}\label{main_009}
\psi(x_t)\equiv x_t\mod{[N,N]};\,\psi(\partial_t(r))\equiv  \partial_t(r)\mod{U_0(N)^i}.
\end{eqnarray}
If $\partial_t(r)$, $\psi(\partial_t(r))\in U_0(N)^i$ then $x_t \partial_t(r)$, $\psi(x_t) \psi(\partial_t(r))\in U_0(N)^{i+1}$.
If $x_t$ --- $\alpha$-generator then from $[N,N]\subseteq U_0(N)^2$, $\partial_t(r)\in U_0(N)^{i-1}$ and (\ref{main_009})
follows that $\psi(x_t) \psi(\partial_t(r))\equiv x_t \partial_t(r)\mod{U_0(N)^{i+1}}$.
Hence, $r -\hat{r}\in U_0(N)^{i+1}$, i.e. $r\in H+N_{1,i+1}$, a contradiction.
So
\begin{eqnarray}\label{tm2_6_gr_end_00}
\partial_z(v)\in \Delta_l,\,z=1,\ldots,a
\end{eqnarray}
as we wished to show. Since $U(H_1)\cap \Delta_l\equiv \Delta_l^\prime\mod{(N_l)_U}$ it follows from (\ref{tm2_6_gr_end_00}) that $\partial_z(v)\in \Delta_l^\prime\mod{(N_l)_U}$. Therefore there exists $v_z\in \Delta_l^\prime$  such that
$\partial_z(v)-v_z\in U(H_1)\cap (N_l)_U$ whence $\partial_z(v)\in U(H_1)(H\cap N_l) + \Delta_l^\prime$ $(z=1,\ldots,a)$.

We have, by the inductive assumption, that $H\cap N_l=H\cap N_{1l}=(H\cap N)_{(l)}\subseteq U_0(H_1)^l$, whence $\partial_z(v)\in \Delta_l^\prime$ $(z=1,\ldots,a)$.
If $k\leqslant l$ then $H_k=H\cap N_{1k} = (H_1)_{(k)}\subseteq U_0(H_1)^k$ hence
$\Delta_l^\prime\subseteq U_0(H_1)^l$, i.e. $\partial_z(v)\in U_0(H_1)^l$, $z=1,\ldots,a$.
Thus $v\in U_0(H_1)^{l+1}$ whence $v\in (H\cap N)_{(l+1)}=H\cap N_{1,l+1}$ and, by induction on $l$,
$H\cap (R+N_{1l}) = H\cap N_{1l}$ for each term $N_{1l}$ of series {\rm (\ref{tm5_0})}.
\end{proof}

From Propositions \ref{tm4}, \ref{tm2}, \ref{tm5}, we obtain immediately Theorem~\ref{tm2_alg}.

\section{A generalized theorem on freedom for relatively free Lie algebras}
Suppose $G$ is a Lie algebra. The elementary transformations
of a matrix over $U(G)$ are defined as follows:
\begin{align}
&\text{Interchange of the columns $i$ and $j$;}\label{gr2_df1}\\
&\text{Interchange of the rows $i$ and $j$;}\label{gr2_df2}\\
&\text{Right-multiply the }i\text{-th row by a non-zero element of }U(G);\label{gr2_df3}\\
&\text{Add to the }j\text{-th row the }i\text{-th row, }\label{gr2_df4}\\
&\text{right-multiplied by a non-zero element of }U(G)\text{, where }i<j.\notag
\end{align}

Let $M=\|m_{kn}\|$ be a $r\times s$ matrix over $U(G)$, $t$ --- the rank of $M$. We denote by $\Phi(M)$ the matrix obtained from $M$ by a chain $\Phi$ of elementary transformations of $M$; $M$ is said to be the lower triangular matrix if $m_{kk}\neq 0\,(k=1,\ldots, t)$, $m_{kn}= 0\,(k>n\,or\, k>t)$.\\
Let $v$ be a $1\times s$ matrix over $U(G)$. We define $\Phi(v)$
by putting for any $\psi$ of $\Phi$ if $\psi$ --- elementary transformation of rows then $\psi(v)=v$.

\begin{lemma}\label{lm4_2_gr_1}
Let $F$ be a  Lie algebra,
$N=N_1 \geqslant \ldots \geqslant N_m\geqslant \ldots$  --- a descending series of ideals
of $F$ such that $[N_p\,,N_q\,]\leqslant N_{p+q}$,
$F/N$ --- a soluble Lie algebra.\\
Let $S=S_\alpha$, $\phi$ --- natural homomorphism $U(F)\to U(F/N_m)$,
$\|a_{kn}\|$ --- $r\times s$ matrix over $U(F/N_m)$, $\phi^\prime$ --- natural homomorphism $U(F/N_m)\to U(F/N)$,  $\phi^\prime(a_{kk})\neq 0$, if $n<k$, then $\phi^\prime(a_{kn})= 0$ $(k=1,\ldots, r)$, $\psi$ --- a mapping of the set $U(F/N_m)$ into
$\{\infty,\,0,\,{\bf N}\}$ such that $\psi(0)=\infty$; if $u\in S\Delta_j\setminus S\Delta_{j+1}\mod{(N_m)_U}$, then $\psi(\phi(u))=j$.
Then matrix $\|a_{kn}\|$ by a finite number of operations {\rm (\ref{gr2_df3})}, {\rm (\ref{gr2_df4})}
can be converted into a matrix $\|b_{kn}\|$ such that $\psi(b_{kk})\leqslant \psi(b_{kn})$; $b_{kk}\neq 0$; if $n<k$, then $b_{kn}= 0$ (k=1,\ldots, r; n=1,\ldots, s).
\end{lemma}
\begin{proof}
Since $F/N_m$ is a soluble Lie algebra, $U(F/N_m)$ is known to satisfy right Ore's condition (see e.g. \cite{Hm}).

If $u\in S\Delta_i\setminus S\Delta_{i+1}\mod{(N_m)_U}$, $v\in S\Delta_j\setminus S\Delta_{j+1}\mod{(N_m)_U}$,
then $uv$ --- an element of $S\Delta_{i+j}\setminus S\Delta_{i+j+1}\mod{(N_m)_U}$.
Hence, $\psi(uv)=\psi(u)+\psi(v)$. It is easily seen that $\psi(u+v)\geqslant\min\{\psi(u),\psi(v)\}$, i.e. $\psi$ --- a valuation on $U(F/N_m)$.

By hypothesis, $\phi^\prime(a_{kk})\neq 0$, if $n<k$, then $\phi^\prime(a_{kn})= 0$.
Hence $\psi(a_{kk})=0$ and if $n<k$ then $\psi(a_{kn})>0$, $k=1,\ldots,r$.

We put $(b_{11},\ldots,b_{1s})=(a_{11},\ldots,a_{1s})$. It is clear that $\psi(b_{11})\leqslant \psi(b_{1j})$, $j=1,\ldots,s$.
Now assume inductively that by a finite number of operations {\rm (\ref{gr2_df3})}, {\rm (\ref{gr2_df4})}
the rows $(a_{k1},\ldots,a_{ks})$ can be converted into $(b_{k1},\ldots,b_{ks})$
such that $b_{kk}\neq 0$; if $n<k$ then $b_{kn}= 0$; $\psi(b_{kk})\leqslant \psi(b_{kn})$, $k=1,\ldots, t-1; n=1,\ldots,s$.

If $a_{t1}=\ldots=a_{t,t-1}=0$ then we put $(b_{t1},\ldots,b_{ts})=(a_{t1},\ldots,a_{ts})$. It is clear that
$b_{tt}\neq 0$; if $n<t$ then $b_{tn}= 0$; $\psi(b_{tt})\leqslant \psi(b_{tj})$, $j=1,\ldots,s$.

Let us now suppose that $a_{t1}=\ldots=a_{t,l-1}=0$, $a_{t,l}\neq 0$, $l\leqslant t-1$. There exist a non-zero elements $\beta_1,\,\beta_2$ of $U(F/N_m)$ such that $b_{ll}\beta_1=-a_{tl}\beta_2$.\\
We put $c_{tn}=b_{ln}\beta_1+a_{tn}\beta_2$, $n=1,\ldots,s$. Then $c_{t1}=\ldots=c_{tl}= 0$.

From $\psi(b_{ll}\beta_1)\leqslant \psi(b_{lj}\beta_1)$ $(j=1,\ldots,s)$, $\psi(b_{ll}\beta_1)=\psi(a_{tl}\beta_2)$
we obtain\\
$\psi(a_{tt}\beta_2)< \psi(a_{tl}\beta_2)\leqslant \psi(b_{lj}\beta_1)$ $(j=1,\ldots,s)$.\\
Since $\psi(a_{tt}\beta_2)< \psi(a_{tj}\beta_2)$ $(j<t)$ and $\psi(a_{tt}\beta_2)\leqslant \psi(a_{tj}\beta_2)$ $(j=1,\ldots,s)$,
we have $\psi(c_{tt})=\psi(a_{tt}\beta_2)< \psi(c_{tj})$ $(j<t)$ and $\psi(c_{tt})\leqslant \psi(c_{tj})$, $j=1,\ldots,s$.

Thus the row $(a_{t1},\ldots,a_{ts})$ can be converted by a finite number of operations {\rm (\ref{gr2_df3})}, {\rm (\ref{gr2_df4})} into $(b_{t1},\ldots,b_{ts})$ such that $b_{tt}\neq 0$; if $n<t$ then $b_{tn}= 0$; $\psi(b_{tt})\leqslant \psi(b_{tj})$ $(j=1,\ldots,s)$ and the proof is complete.
\end{proof}

Let $\|a_{kn}\|$ be a $r\times s$ matrix over $U(F/N_m)$, $t$ --- the rank of $\|a_{kn}\|$.
It is not hard to verify that $\|a_{kn}\|$ can be converted by a finite number of operations
{\rm (\ref{gr2_df1})}, {\rm (\ref{gr2_df2})},
{\rm (\ref{gr2_df3})}, {\rm (\ref{gr2_df4})}
into a lower triangular form $\|b_{kn}\|$ such that $\psi(b_{kk})\leqslant \psi(b_{kn})\,(k=1,\ldots, t; n=1,\ldots, s)$.

\begin{lemma}\label{lm4_2_gr_2}
Suppose $G$ is a soluble Lie algebra; $M$ --- $r\times s$ matrix over $U(G)$; $\alpha_i$ --- $i$-th row of $M$;
$\alpha$ --- a right-linear combination of rows $\alpha_1,\ldots,\alpha_r$.
Then for each elementary transformation $\psi$ of $M$ there exists non-zero $d_\psi\in U(G)$ such that $\psi(\alpha) d_\psi$ is a right-linear combination of rows of $\psi(M)$.
\end{lemma}
\begin{proof}
If $\psi$ is one of the operations  (\ref{gr2_df1}),
(\ref{gr2_df2}), (\ref{gr2_df4}) the result is obvious.
Let $\psi(M)$ be a matrix obtained by operation (\ref{gr2_df3}); $a$ ---
non-zero element of $U(G)$; $\alpha_1,\ldots,\alpha_i a,\ldots,\alpha_r$ --- the rows of
$\psi(M)$.
By assumption, there exist an elements $b_1,\ldots,b_r$ of $U(G)$ such that
$\alpha_1b_1+\ldots+\alpha_ib_i+\ldots+\alpha_rb_r=\alpha$. It is no restriction to
assume that $b_i\neq 0$, otherwise the result is obvious.
When $G$ is a soluble Lie algebra, $U(G)$ is known to satisfy right Ore's condition.
Thus there exist non-zero $c$, $d_\psi\in U(G)$ such that $ac=b_id_\psi$ and we have $\alpha_1b_1d_\psi+\ldots+\alpha_i ac+\ldots+\alpha_rb_rd_\psi=\alpha\, d_\psi$.
\end{proof}

The proof of Theorem \ref{tm3_alg}. We may clearly assume that $n-m>1$. Let $\phi_k$  be the natural homomorphism
$U(F)\to U(F/(R+N_{k,m_k+1}))$, $\phi_k^\prime$ --- natural homomorphism $U(F/(R+N_{k,m_k+1}))\to U(F/(R+N_{k1}))$,
$A=\|a_{rt}\|$ --- a matrix over $U(F)$. We denote by $A^{\phi_k}$ the matrix $\|\phi_k(a_{rt})\|$
and denote by $(A^{\phi_k})^{\phi_k^\prime}$ the matrix $\|\phi_k^\prime(\phi_k(a_{rt}))\|$.

We denote by $N$ the ideal $R+N_{k1}$ and denote by $N_p$ the ideal $R+N_{kp}$, $p\in {\bf N}$.
It is not hard to verify that $[N_p\,,N_q\,]\leqslant N_{p+q}$.
We define a valuations $\psi_0$ on $U(F/N_{11})$ and $\psi_k$ on $U(F/N_{m_k+1})$ as follows:
$\psi_k(0)=\psi_0(0)=\infty$; if $u\neq 0$ then $\psi_0(u)=0$; if $u\in S\Delta_j\setminus S\Delta_{j+1}\mod{(N_{m_k+1})_U}$ then $\psi_k(\phi_k(u))=j$.

Let $D_1,\ldots,D_n$ be the Fox derivatives of $U(F)$.
We denote by $m_{ij}$ the elements $D_j(r_i)\,(i=1,\ldots, m;\,j=1,\ldots, n)$, by $M$ --- the matrix $\|m_{ij}\|$, by $t_k$ --- the rank of $M^{\phi_k}$.

If $m_{ij}\equiv 0 \mod{(R+N_{s,m_s+1})_U}$ for any $m_{ij}$ then $r_i\in (R+N_{s,m_s+1})_{(2)}\,(i=1,\ldots, m; j=1,\ldots, n)$ \cite{Hm}. From $r_i\in (R+N_{s,m_s+1})_{(2)}$ (i=1,\ldots, m) follows that $R\subseteq (R+N_{s,m_s+1})_{(2)}$, whence
$R\subseteq [R,R]+N_{s,m_s+1}$.
Since $F/N_{s,m_s+1}$ --- a solvable algebra we have $R\subseteq N_{s,m_s+1}$.
Thus $H\cap (R+N_{kl}) = H\cap N_{kl}$ for each term $N_{kl}$ of series (\ref{end_algr_3}), where $H=F$.

Further we may assume that there exists $K\in \{0,\ldots,s\}$ such that $t_K>0$ and if $i<K$ then $t_i=0$.
Let $\Phi_K$ be a chain of elementary transformations of $M$ such that $M_K=\|m^{(K)}_{ij}\|=(\Phi_K(M))^{\phi_K}$ be a lower triangular $m\times n$ matrix and $\psi_K(m^{(K)}_{ii})\leqslant \psi_K(m^{(K)}_{ij})\,(i=1,\ldots, t_K;\,j=1,\ldots, n)$.

Let $K<s$. Now assume inductively that for some $k$ $(K\leqslant k< s)$ we have $\Phi_k$ --- a chain of elementary transformations of $M$;
$M_k=\|m^{(k)}_{ij}\|=(\Phi_k(M))^{\phi_k}$ --- a lower triangular $m\times n$ matrix;
$\psi_k(m^{(k)}_{ii})\leqslant \psi_k(m^{(k)}_{ij})\,(i=1,\ldots, t_k;\,j=1,\ldots, n)$.

We denote $\Phi_k(M)$ by $M_{k+1,1}$. We have $(M_{k+1,1}^{\phi_{k+1}})^{\phi^\prime_{k+1,1}}=M_k$. Then Lemma~\ref{lm4_2_gr_1} tells us that there exists a chain $\Phi_{k+1,1}$ of elementary transformations (\ref{gr2_df3}) (where $i\leqslant t_k$), (\ref{gr2_df4}) (where $i<j\leqslant t_k$) of $M_{k+1,1}$ such that $(\Phi_{k+1,1}(M_{k+1,1}))^{\phi_{k+1}}$ be $m\times n$ matrix $\|b_{ij}\|$ with $b_{ii}\neq 0$; $b_{ij}=0\,(j< i)$; $\psi_{k+1}(b_{ii})\leqslant \psi_{k+1}(b_{ij})\,(i=1,\ldots, t_k;\,j=1,\ldots, n)$.

We denote $\Phi_{k+1,1}(M_{k+1,1})$ by $M_{k+1,2}$
and denote by $\Phi_{k+1,2}$ a chain of operations (\ref{gr2_df3}) (where $i> t_k$), (\ref{gr2_df4}) (where $i\leqslant t_k\text{ and }j>t_k$) such that $(\Phi_{k+1,2}(M_{k+1,2}))^{\phi_{k+1}}$ be a matrix $\|c_{ij}\|$ with $c_{ij}=b_{ij}\,(i=1,\ldots, t_k;\,j=1,\ldots, n)$ and $c_{ij}=0\,(j\leqslant t_k<i)$.

We denote $\Phi_{k+1,2}(M_{k+1,2})$ by $M_{k+1,3}$
and denote by $\Phi_{k+1,3}$ a chain of operations (\ref{gr2_df1}), (\ref{gr2_df2}), (\ref{gr2_df3}), (\ref{gr2_df4}) (where $i,\,j> t_k$) of $M_{k+1,3}$ such that $(\Phi_{k+1,3}(M_{k+1,3}))^{\phi_{k+1}}$ be a lower triangular $m\times n$ matrix $\|m^{(k+1)}_{ij}\|$ and $\psi_{k+1}(m^{(k+1)}_{ii})\leqslant \psi_{k+1}(m^{(k+1)}_{ij})\,(i=1,\ldots, t_{k+1};\,j=1,\ldots, n)$.

We denote the matrix $\|m^{(k+1)}_{ij}\|$ by $M_{k+1}$, the sequence $\Phi_k$, $\Phi_{k+1,1}$,
$\Phi_{k+1,2}$, $\Phi_{k+1,3}$ by $\Phi_{k+1}$.
Thus by induction on $k$ we have $M_k$, $\Phi_k$ for any $k\in \{K,\ldots,s\}$.
If $K=s$ the truth of the statement is obvious.

Let $I_s=\{i_1,\ldots,i_{t_s}\}$ be the subset of $\{1,\ldots,n\}$ such that
if $m_{i_j}$ --- $i_j$-th column of $M$ then $\Phi_s(m_{i_j})$ --- $j$-th column of $\Phi_s(M)$;
$\{j_1,\ldots,j_p\}= \{1,\ldots,n\}\setminus I_s$;
$H$ --- the free Lie algebra with the free set  $\{y_{j_1},\ldots,y_{j_p}\}$ of generators.
Since $t_s\leqslant m$ it follows that $p\geqslant n-m$.

It is clear that $H\cap (R+N_{11}) = H\cap N_{11}$.
Now assume inductively that $H\cap (R+N_{ij}) = H\cap N_{ij}\,(i=1,\ldots, k;\,j=1,\ldots, l;\,l\leqslant m_k)$.
We need to show that $H\cap (R+N_{k,l+1})=H\cap N_{k,l+1}$.

Let $\{x_{kz} \mid  z\in I\}$ be a free set of generators of $N$;
$\{\partial_{kz} \mid  z\in I\}$ ---  the Fox derivatives of $U(N)$;
$v\in H\cap N_{m_k+1}$. From $N_{m_k+1}=R+N_{k,m_k+1}=R+N_{(m_k+1)}$ follows the existence of an element $u$ of $N_{(m_k+1)}$ such
that $v-u\in \mbox{\rm id}_F (r_1,\ldots,r_m)$.
From (\ref{four}) follows the existence of an elements $C_1,\ldots,C_m$ of $U(F)$ such
that $D_j(v-u)\equiv \sum_{i={1}}^{m} D_j(r_i)C_i\mod{R_U},\,j= 1,\ldots,n$.
Hence
\begin{eqnarray}\label{tm2_3_gr_end}
D_j(v)\equiv \sum_{i={1}}^{m} D_j(r_i)C_i+D_j(u)\mod{(N_{m_k+1})_U},\,j= 1,\ldots,n.
\end{eqnarray}
Since $u\in N_{(m_k+1)}$ it follows that $\sum_{z\in I} D_j(x_{kz})\partial_{kz}(u)\in U(F)U_0(N)^{m_k}$, i.e.
\begin{eqnarray}\label{tm2_3_gr_end_1}
D_j(u)\in U(F)U_0(N)^{m_k},\,j=1,\ldots,n.
\end{eqnarray}
We now show that from $v\in H\cap N_{m_k+1}$ follows
\begin{eqnarray}\label{tm2_6_gr_end}
D_j(v)\in U(F)\Delta_{m_k}\mod{(N_{m_k+1})_U},\,j=1,\ldots,n.
\end{eqnarray}

If $t_k=0$ then $D_j(r_i)\equiv 0\mod{(N_{m_k+1})_U}\,(i=1,\ldots, m;\,j=1,\ldots, n)$, hence from
(\ref{tm2_3_gr_end}), (\ref{tm2_3_gr_end_1}) follows (\ref{tm2_6_gr_end}).
Let $t_k>0$, $V=(D_1(v)-D_1(u),\ldots,D_n(v)-D_n(u))$.
From Lemma~\ref{lm4_2_gr_2} and (\ref{tm2_3_gr_end}) we obtain the existence of an element $d$ of $U(F)$ such that $\phi_k(d)\neq 0$ and $(\Phi_k(V d))^{\phi_k}$ --- a right-linear combination of non-zero rows of $M_k$.

Since $v\in H$ it follows that $D_j(v)=0$
$(j=i_1,\ldots,i_{t_s})$ and therefore an element of $\Phi_k(V d)$ in the $i$-th coordinate belongs to
$U(F)\Delta_{m_k}d\, (i=1,\ldots, t_k)$.
Since $\psi_k(m^{(k)}_{ii})\leqslant \psi_k(m^{(k)}_{ij})\,(i=1,\ldots, t_k;\,j=1,\ldots, n)$ it is not hard to verify that all elements of $(\Phi_k(V d))^{\phi_k}$ belong to $\phi_k(U(F)\Delta_{m_k}d)$.

Thus $D_j(v)-D_j(u)\in U(F)\Delta_{m_k}\mod{(N_{m_k+1})_U}\,(j=1,\ldots,n)$, whence $D_j(v)\in U(F)\Delta_{m_k}\mod{(N_{m_k+1})_U}\,(j=1,\ldots,n)$ as we wished to show.

Let $v\in H\cap N_{l+1}$. Consider case $l=1$.
We have $H\cap N=H\cap N_{k1}$ and $v\in H\cap N_2$.
From $N_2=R+N_{k2}$ follows the existence of an element $u$ of $N_{k2}$ such
that $v-u\in \mbox{\rm id}_F (r_1,\ldots,r_m)$.
Then from (\ref{four}) follows the existence of an elements $B_1,\ldots,B_m$ of $U(F)$ such
that
\begin{eqnarray}\label{tm2_6_gr_ends}
D_j(v-u)\equiv \sum_{i={1}}^{m} D_j(r_i)B_i\mod{N_U},\,j= 1,\ldots,n.
\end{eqnarray}
Let $V=(D_1(v)-D_1(u),\ldots,D_n(v)-D_n(u))$.
Since $v\in H$ it follows that $D_j(v)=0$
$(j=i_1,\ldots,i_{t_s})$. From $u\in N_{k2}$ we have $\phi_{k-1}((D_j(u))=0$ $(j=1,\ldots,n)$.
From Lemma~\ref{lm4_2_gr_2} and (\ref{tm2_6_gr_ends}) we obtain the existence of an element $d$ of $U(F)$ such that
$\phi_{k-1}(d)\neq 0$ and $(\Phi_{k-1}(V d))^{\phi_{k-1}}$ --- a right-linear combination of rows of a lower triangular matrix $M_{k-1}$. Thus $(\Phi_{k-1}(V d))^{\phi_{k-1}}$ --- trivial row whence
$D_j(v)\equiv 0 \mod{(N_{k1})_U}$ $(j=1,\ldots,n)$ and $v\in N_{k2}$ \cite{Hm}.

Consider case $l>1$. We begin by noting that if $x$ --- an element of a free set of generators of $H\cap N$ then there are an elements $q_x\in \{j_1,\ldots,j_p\}$, $n_x\in {\bf N}$ such that $D_{q_x}(x) \equiv \sum_{i=1}^{n_x} \lambda_i f_i\mod{N_U}$ where $0\neq \lambda_i\in P$, $f_i \in S_\alpha$. It follows because if $D_j(x)\equiv 0\mod{N_U}$, $j= 1,\ldots,n$ then $x\in N_{(2)}$ \cite{Hm} and the above arguments (case $l=1$) give $x\in H\cap N_{k2}=(H\cap N)_{(2)}$, a contradiction.

Let $z_1,\ldots,z_a\,(z_i\neq z_j\, for\, i\neq j) $ be an elements of a free set of generators of $H\cap N$ such that $v$ belong to $(z_1,\ldots,z_a)$.
It is not hard to verify that $\{z_1,\ldots,z_a\}$ by a finite number of elementary transformations
can be converted into a $\{x_1,\ldots,x_a\}$ such that:
$v$ belong to $(x_1,\ldots,x_a)$;
$D_{q_i}(x_i) \equiv \lambda_i M_i+V_{ii}\mod{N_U}$, $D_{q_i}(x_j) \equiv V_{ij}\mod{N_U}\,(
for\, i\neq j)$ where $q_i\in \{j_1,\ldots,j_p\}$, $0\neq \lambda_i\in P$, $M_i \in S_\alpha$,
$V_{ij}$ --- linear combination of monomials of the form (\ref{alg2_1}) and this monomials are not in $M_i U(N)\,(j=1,\ldots,a)$.
We denote by $\partial_1,\ldots,\partial_a$ the Fox derivatives of $(x_1,\ldots,x_a)$.

Since a free set of generators of a free Lie algebra $H$ contains more than one element and
$H\cap N$ --- an ideal of $H$ it follows that a free set of generators of  a free Lie algebra $H\cap N$ contains more than one element.
So we may assume that $a>1$.

From $N_l\subseteq N_2$ and (\ref{fm2_pr_gr}) we have $U(H_1)\cap \Delta_2\equiv \Delta_2^\prime\mod{(N_2)_U}$.\\
Let $y\in U(H_1)\cap \Delta_2\mod{(N_{m_k+1})_U}$.
Then $y\in U(H_1)\cap \Delta_2\mod{(N_2)_U}$.\\
Hence
\begin{eqnarray}\label{z}
y\in\Delta_2^\prime\mod{(N_2)_U}.
\end{eqnarray}
From $H\cap N_2=H\cap N_{k2}=(H\cap N_{k1})_{(2)}\subseteq U_0(H_1)^2$ follows that $\Delta_2^\prime = U_0(H_1)^2$.
From (\ref {z}) follows that $y\in U_0(H_1)^2+U(H_1)\cap (N_2)_U$ whence $y\in U_0(H_1)^2+U(H_1)(H\cap N_2)$ and $y\in U_0(H_1)^2$.
We now show that
\begin{eqnarray}\label{tm2_6_gr_end_0}
\partial_z(v)\in \Delta_l\mod{(N_{m_k+1})_U},~z=1,\ldots,a.
\end{eqnarray}

Suppose that there exists $i$ such that $\partial_i(v)\in \Delta_{l-1}\setminus \Delta_l\mod{(N_{m_k+1})_U}$.
We denote by $\bar{v}$ the element $[\ldots [v,x_t],\ldots ,x_t]$ $(t\neq i)$ such that $\bar{v}\in H\cap N_{m_k+1}$ and $\partial_i(\bar{v})=\partial_i(v)x_t^{m_k -l}$.
Since  $x_t$ --- an element of a free set of generators of $H_1$ it follows that $x_t\notin U_0(H_1)^2$.
Hence $x_t\notin \Delta_2\mod{(N_{m_k+1})_U}$.
From $x_t\not\in \Delta_2\mod{(N_{m_k+1})_U}$ and
$\partial_i(\bar{v})=\partial_i(v)x_t^{m_k -l}$ follows that
$\partial_i(\bar{v})\in \Delta_{m_k -1}\setminus \Delta_{m_k}\mod{(N_{m_k+1})_U}$.

By (\ref{alg1_dcf}), $D_{q_i}(\bar{v})=\sum_{z={1}}^{a} D_{q_i}(x_z)\partial_z(\bar{v})$, whence
\begin{eqnarray*}
D_{q_i}(\bar{v}) = \lambda_i\cdot M_i\partial_i(\bar{v})+V,
\end{eqnarray*}
where $0\neq \lambda_i \in P$, $V$ --- linear combination of monomials of the form (\ref{alg2_1})
and this monomials are not in $M_i U(N)$.

Thus we have $D_{q_i}(\bar{v})\not\in U(F)\Delta_{m_k}\mod{(N_{m_k+1})_U}$
which contradicts (\ref{tm2_6_gr_end}).
This proves (\ref{tm2_6_gr_end_0}).

Since $U(H_1)\cap \Delta_l\equiv \Delta_l^\prime\mod{(N_l)_U}$ and $(N_{m_k+1})_U\subseteq(N_l)_U$ it follows from (\ref{tm2_6_gr_end_0}) that $\partial_z(v)\in \Delta_l^\prime\mod{(N_l)_U}$.
Therefore there exists $v_z\in \Delta_l^\prime$  such that
$\partial_z(v)-v_z\in U(H_1)\cap (N_l)_U$ whence $\partial_z(v)\in U(H_1)(H\cap N_l) + \Delta_l^\prime$ $(z=1,\ldots,a)$.

We have, by the inductive assumption, that $H\cap N_l=H\cap N_{kl}=(H\cap N_{k1})_{(l)}\subseteq U_0(H_1)^l$, whence $\partial_z(v)\in \Delta_l^\prime$.
If $t\leqslant l$ then $H_t=H\cap N_{kt} = (H_1)_{(t)}\subseteq U_0(H_1)^t$ hence
$\Delta_l^\prime\subseteq U_0(H_1)^l$, i.e. $\partial_z(v)\in U_0(H_1)^l$, $z=1,\ldots,a$.
Thus $v\in U_0(H_1)^{l+1}$ whence $v\in (H\cap N_{k1})_{(l+1)} = H\cap N_{k,l+1}$ and, by induction on $l$,
$H\cap (R+N_{kl}) = H\cap N_{kl}$ for each term $N_{kl}$ of series {\rm (\ref{end_algr_3})}.



\begin{thebibliography}{99}

\bibitem{Sh} A.I.Shirshov, Some algorithmic problems for Lie algebras, Sibirsk. Mat. \v{Z}., 3, N~2 (1962), 292--296.

\bibitem{Tl} V.V.Talapov, On polynilpotent Lie algebras defined by one defining relation, Sibirsk. Mat. \v{Z}., 23, N~5 (1982), 192--204.

\bibitem{Hm} O.G.Kharlampovich, Lyndon condition for solvable Lie algebras, Izv. Vyssh. Uchebn. Zaved. Mat., N~9 (1984), 50-59.

\bibitem{Sh2} A.I.Shirshov, Subalgebras of free Lie algebras, Mat. Sb., 33, N~2 (1953), 441--452.

\end{thebibliography}
\end{document}